\documentclass{amsart}
\usepackage{amssymb}
\usepackage{amscd}
\usepackage{epsfig}
\usepackage{float}

\newtheorem{thm}{Theorem}[section]
\newtheorem{lem}[thm]{Lemma}
\newtheorem{defn}[thm]{Definition}
\newtheorem{cor}[thm]{Corollary}
\newenvironment{pf}{\noindent \textbf{Proof.}}{\qed\vspace{.5\baselineskip}}

\newcommand{\field}[1]{\mathbb{Z}_{#1}}
\newcommand{\norm}[1]{x(#1)}
\newcommand{\intersecto}{\triangle}
\newcommand{\intersect}[2]{\triangle(#1,#2)}
\newcommand{\lk}[2]{\operatorname{lk}(#1,#2)}
\newcommand{\trace}{\operatorname{trace}}
\newcommand{\image}{\operatorname{image}}

\newcommand{\col}[2]{(\begin{smallmatrix}#1 \\ #2 \end{smallmatrix})}

\def\R{\mathbb{R}}
\def\Z{\mathbb{Z}}
\def\Q{\mathbb{Q}}
\def\d{\partial}
\def\Ft{{\tilde F}}
\def\ft{{\tilde f}}
\def\Kt{{\tilde K}}

\def\Mt{{\tilde M}}
\def\St{{\tilde S}}
\def\Tt{{\tilde T}}

\def\xt{{\tilde x}}
\def\yt{{\tilde y}}
\def\alphat{{\tilde \alpha}}

\newsavebox{\ptorus}
\savebox{\ptorus}{\put(0,0){\framebox(20,20){}}
\put(10,10){\circle{5}}}

\begin{document}

\title{Virtually Embedded Boundary Slopes}
\author{Joseph Maher}
\date{October 15, 1997}

\begin{abstract}
We show that for certain hyperbolic manifolds all boundary slopes are 
slopes of \text{$\pi_{1}$-injective} immersed surfaces, covered by incompressible embeddings in 
some finite cover. The manifolds include hyperbolic punctured torus 
bundles and hyperbolic two-bridge knots.


Subject code: 57M10, 57M25

\end{abstract}

\maketitle

\section{Introduction}

Let $M$ be a manifold with boundary consisting of a single torus, 
$T$. A \emph{slope}, $\alpha$, is an essential simple closed 
curve in $T$. The slope is an \emph{embedded boundary slope} if there 
is an embedded surface in $M$, whose boundary consists of loops in 
$T$ parallel to $\alpha$. The surface must be compact and orientable. 
It must also be properly embedded, $\pi_{1}$-injective, and not 
properly homotopic rel boundary to any part of the boundary of $M$.

We can define boundary slopes for manifolds with more than one torus boundary component in the same way. Then a properly embedded surface will have a boundary slope defined for each boundary component of the manifold it intersects. If the surface does not intersect a particular boundary component, then the slope of the surface on that component is not defined.

If we only require the surface to be an immersion, though still 
embedded in a neighbourhood of $T$, then the slope is an 
\emph{immersed boundary slope}. If the immersion is covered by an 
embedding into some finite cover of $M$, then the slope is a 
\emph{virtually embedded boundary slope}.

It has been shown that a knot can have only finitely many embedded 
boundary slopes \cite{Hatcher82}, and many examples of such surfaces 
have been constructed, for example \cite{HO}. Examples of immersed 
boundary slopes have been found in the figure eight knot 
\cite{Hempel86,Hempel87}. Also \cite{BC} show that for this knot 
every 
slope with even numerator is a virtually embedded boundary slope. It 
has been shown that a punctured torus bundle may have infinitely many 
virtually embedded boundary slopes \cite{Baker}, and that there is a  
manifold with every slope an immersed boundary slope \cite{Oertel}. 

Note that there will always be one slope on $\d M$ which is null 
homologous in $H_{1}(M;\R)$, label this slope $l$. Choose a slope 
$m$, such that $m \cap l$ is a single point. These two curves form a 
basis for $H_{1}(\d M;\R)$. If $M$ is a knot space, choose $l$ and 
$m$ to be the longitude and meridian of the knot. We will label the 
slope $a/b$ if the slope is made up of $a$ meridians and $b$ 
longitudes.

In this paper we show:

\begin{thm}{} \label{main}
If $M$ has hyperbolic interior, and a finite cover $\Mt$ such that:
\begin{list} {\emph{\arabic{enumi}.}}{\usecounter{enumi}}
\item $\Mt$ has at least three boundary components.
\item There is a boundary torus $\Tt$ of $\Mt$ which is a one-fold 
covering of $\d M$.
\item The projection $\rho :\ker i_{*} \rightarrow H_{1}(\Tt;\R)$ is 
onto, where $i_{*}:H_{1}(\d \Mt;\R)\rightarrow H_{1}(\Mt;\R)$ is the 
map induced by inclusion, and $\rho$ is the vector space projection $\rho:H_{1}(\d 
\Mt;\R)\rightarrow H_{1}(\Tt;\R)$.
\end{list}
Then every slope of $\d M$ is a virtually embedded boundary slope. 
\end{thm}

Note that in condition 3, the map $\rho$ is projection onto $H_{1}(\Tt;\R)$, which is a {\em subspace} of $H_{1}(\d \Mt;\R)$, and this map has nothing to do with the covering projection $p:\Mt \rightarrow M$, and in fact is not induced by any continuous map between manifolds.

In \cite{BC} it is shown that if there is a 
surface with boundary slope $a/b$ on $\Tt$, which is not the fiber of a 
fibration, then $a/b$ is a virtually embedded boundary slope of $\Tt$. 
As $\Tt$ is a one-fold covering of $\d M$, these surfaces project down 
to 
virtually embedded surfaces in $M$ with the same boundary slope. We 
then use the Thurston norm to show that there are surfaces of every 
boundary slope on $\Tt$, which are not fibers of fibrations.

By constructing particular covers we then show:

\begin{cor}{} \label{torus}
For hyperbolic punctured torus bundles, every boundary slope is a 
virtually embedded boundary slope.
\end{cor}

\begin{cor} \label{s2s1}
If $K$ is a knot in $S^{2} \times S^{1}$ such that $M = S^{2} \times 
S^{1} - K$ is hyperbolic, and the algebraic intersection number of 
$K$ with $S^{2} \times \{point \}$ is at least three, then every 
slope of $M$ is an embedded boundary slope.
\end{cor}

\begin{cor} \label{twobridge}
For hyperbolic two bridge knots, every boundary slope is a virtually 
embedded boundary slope.
\end{cor}

\section{General Discussion.}

\subsection{Proof of Theorem \ref{main}.}

Let $M$ be a hyperbolic $3$-manifold with boundary consisting of a 
single torus, and let $\Mt$ be a finite cover of $M$. Note that $\Mt$ 
will also be hyperbolic and hence irreducible, and that $\d \Mt$ will 
be incompressible in $\Mt$.

The main result we will use is: 

\begin{thm}{\cite[Theorem 1.4]{BC}} \label{bc}
Let $\Mt$ be a compact, connected, orientable, atoroidal and 
irreducible 3-manifold, with boundary a finite number of tori. 
Suppose that $S$ is a connected, non-separating, orientable, 
incompressible surface properly embedded in $\Mt$, which is not the 
fiber of a fibration of $\Mt$. Also suppose that $\d S$ contains some 
components with slope $\alpha$, on a torus, $\Tt$, in the boundary of 
$\Mt$. Then $\alpha$ is a virtually embedded boundary slope.

Moreover there is a finite cyclic cover $\Mt_{2}$, and a 
compact, connected, orientable, incompressible, 
boundary-incompressible, surface $F$, properly embedded in $\Mt_{2}$. 
The 
boundary of $F$ consists of a non-empty set of essential, parallel 
curves lying on some component $U$ of $\d \Mt_{2}$ which covers 
$\Tt$. 
Also the covering map is an immersion on $F$, which is an embedding 
in a neighbourhood of the boundary, and the boundary of F is mapped 
to 
loops parallel to $\alpha$.
\end{thm}

This shows that it suffices to find a connected surface $S$, in some 
finite 
cover $\Mt$ of $M$, with the following properties:

\begin{list}{(\roman{enumi})}{\usecounter{enumi}}
\item $S$ has the required boundary slope on the boundary torus $\Tt$ 
of 
$\Mt$.
\item $S$ is not the fiber of a fibration.
\end{list}

We will need to know when a surface is not a fiber of a fibration. 
For this we will use the following results about the Thurston norm: 

\begin{defn}{\cite[Section 1, pp 103-105]{Thurston}} If $S$ is 
a connected surface, let $\chi_{-}(S)=max\{0, -\chi(S)\}$, where 
$\chi(S)$ is the 
Euler characteristic of $S$. If $S$ has connected components $S_{1}, 
\ldots, S_{k}$, define 
$\chi_{-}(S)=\chi_{-}(S_{1})+\ldots+\chi_{-}(S_{k})$.

Define the \emph{\textbf{Thurston norm}} $\norm{.}$ on $H_{2}(M,\d 
M;\Z)$ by \[\norm{s}=min\{\chi_{-}(S) \mid S \text{ is an embedded 
surface representing } s \in H_{2}(M,\d M;\Z) \}\]

This extends to a function on $H_{2}(M,\d M;\Q)$ by linearity, and 
then to a function on $H_{2}(M,\d M;\R)$ by continuity.
\end{defn}

\begin{thm}{\cite[Theorem 1, p $100$]{Thurston}} The function 
$\norm{.}$ defined on $H_{2}(M,\d M;\R)$ is convex and linear on rays 
through the origin. If every embedded surface representing a non-zero 
element of $H_{2}(M,\d M;\R)$ has negative Euler characteristic, then 
$\norm{.}$ is a norm. In general $\norm{.}$ is a pseudonorm vanishing 
on precisely the subspace spanned by embedded surfaces of 
non-negative Euler characteristic. \end{thm}

\begin{thm}{\cite[Theorem 2, p $106$]{Thurston}} \label{polyhedron} 
When $\norm{.}$ 
is a norm, the unit ball of $H_{2}(M,\d M;\R)$ is a polyhedron 
defined by linear inequalities with integer coefficients, with 
respect to a basis of primitive elements of $H_{2}(M,\d M;\Z)$.
\end{thm} 

\begin{thm}{\cite[Theorem 3, p $113$]{Thurston}} If the fiber of 
$M$ is a surface with negative Euler characteristic, then the ray 
determined by the homology class of any fiber passes through the 
interior of a top dimensional face of the unit sphere.
\end{thm} 

First we show that $\norm{.}$ is a norm on $H_{2}(\Mt,\d \Mt;\R)$.

\begin{lem}If $\Mt$ satisfies the conditions given above, then the 
Thurston norm $\norm{.}$ is a norm on $H_{2}(\Mt,\d \Mt;\R)$.
\end{lem}

\begin{pf} It suffices to show that if $\norm{s}=0$, then $s$ is 
zero in $H_{2}(\Mt,\d \Mt;\R)$.

If $\norm{s}=0$, then $s$ is represented by an embedded surface $S$, 
whose connected components must all be discs, annuli, spheres, or 
tori.

An embedded sphere in $\Mt$ bounds a ball, as $\Mt$ is irreducible, 
so all embedded spheres represent trivial homology in $H_{2}(\Mt,\d 
\Mt;\R)$.

A properly embedded disc $D$ must bound a disc $D'$ in $\d \Mt$, as 
the boundary of $\Mt$ is incompressible. These two discs bound a 
ball, as $\Mt$ is irreducible, so $D$ represents a trivial homology 
class in $H_{2}(\Mt,\d \Mt;\R)$ as well.

Suppose $S$ has a torus component $Y$. If $Y$ is $\pi_{1}$-injective 
then $Y$ must be parallel to a boundary component of $\Mt$, as $\Mt$ 
is hyperbolic. Therefore $Y$ represents trivial homology.

If $T$ is not $\pi_{1}$-injective, then there is a compressing disc 
$D$ for $Y$ in $\Mt$. Surger $Y$ along $D$ to produce a sphere, which 
bounds a ball in $\Mt$ by irreducibility. So $Y$ is a boundary, as 
$Y$ bounds the union of 
this ball and a regular neighbourhood of $D$, so $Y$ must be trivial 
in 
$H_{2}(\Mt,\d \Mt;\R)$.

Suppose a component of $S$ is an annulus, $A$. If a component of $\d 
A$ bounds a disc $D$ in $\d \Mt$, then $D$ can be pushed off the 
boundary of $\Mt$ so that $A \cup D$ is a properly embedded  disc in 
$\Mt$. This bounds a region $R$ in $(\Mt, \d \Mt)$. So $A$ bounds the 
union of $R$ and a regular neighbourhood of $D$ in $(\Mt, \d \Mt)$, 
so 
$A$ is trivial in $H_{2}(\Mt,\d \Mt;\R)$. 

If both components of $\d A$ are contained in the same boundary 
component $Y$ of $\Mt$,  then the annulus $A$ forms a homotopy in 
$\Mt$ between them. As $Y$ is $\pi_{1}$-injective, and neither curve 
is trivial, they must bound an annulus $A'$ in $Y$. Form a new torus 
$Y'$ from $Y$ by replacing $A'$ with $A$. The embedded torus $Y'$ in 
$\Mt$ is trivial, so $A$ must be trivial in $H_{2}(\Mt,\d \Mt;\R)$ as 
well. 

Suppose the components of $\d A$ are contained in different boundary 
components, $Y_{1}$ and $Y_{2}$, of $\Mt$. Let $N$ be a regular 
neighbourhood of $Y_{1} \cup A \cup Y_{2}$. As neither component of 
$\d A$ bounds a disc in $\d \Mt$, $\d N$ is an embedded torus in 
$\Mt$, so it must be parallel to a boundary component. This means the 
manifold must be a solid torus with two parallel cores drilled out, 
as in Figure \ref{solidtorus}.

\begin{figure}[H]
\begin{center}
\epsfig{file=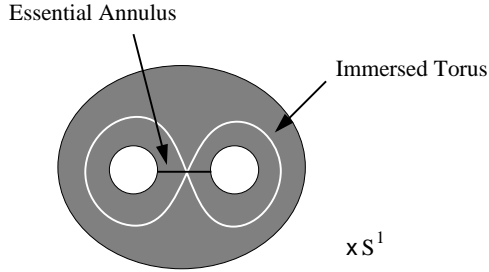, height=100pt}
\caption{Cross section of the manifold} \label{solidtorus}
\end{center}
\end{figure}

However this manifold is not hyperbolic as it contains an immersed 
torus which is not boundary parallel. So in fact no such essential 
annuli can exist.

This shows that if $\norm{s}=0$, then $s=0$ in $H_{2}(\Mt,\d 
\Mt;\R)$, so $\norm{.}$ is a norm.
\end{pf}

Suppose $\Mt$ and $\Tt$ are as in Theorem \ref{main}. We can use the 
Thurston norm to find surfaces which are not fibers, 
using the relative homology exact sequence, as follows:

\begin{equation*}
\begin{CD}
H_{2}(\Mt,\d \Mt;\R) @>\d_{*}>> H_{1}(\d \Mt;\R) @>i_{*}>> 
H_{1}(\Mt;\R) \\
 & & @V \rho VV \\
 & & H_{1}(\Tt;\R) \cong \R^{2} \\
\end{CD}
\end{equation*}

Here $\rho$ is the projection $\rho:H_{1}(\d \Mt;\R) \rightarrow 
H_{1}(\Tt;\R)$, where $\Tt$ is the boundary torus which is a one-fold 
covering of $\d M$.

A boundary slope $a/b$ defines a $1$-dimensional subspace, which we 
will call a  line, in $H_{1}(\Tt;\R)$. We say a line has 
\emph{rational slope} if it contains a non-zero integer homology 
class. We would like to show that the pre-image of this subspace in 
$H_{2}(\Mt,\d \Mt; \R)$ contains a line with rational slope, which 
does not pass through the interior of a top dimensional face of the 
unit ball.

To do this, we will use the following result:

\begin{lem}{} \label{onto}
Suppose $\Mt$ is hyperbolic, with at least three boundary components. 
If the linear map $\phi=\rho \circ \d_{*}:H_{2}(\Mt, \d 
\Mt;\R)\rightarrow 
H_{1}(\Tt;\R)$ is onto, then for any line 
$L$ in $H_{1}(\Tt;\R)$, there is a line $V$ in $\phi^{-1}(L)$, such 
that $\phi(V)=L$, and $V$ does 
not pass through the interior of a top dimensional face of the unit 
ball under the Thurston norm.

Furthermore, if the line $L$ has rational slope, then $V$  can be 
chosen to have rational slope as well. 
\end{lem}

\begin{pf}
By Poincar\'e duality, and the relative homology long exact sequence, we know that the dimension of the kernel of the inclusion map $i_{*}:H_{1}(\d \Mt;\R) \rightarrow H_{1}(\Mt;R)$ is half the dimension of $H_{1}(\d \Mt;\R)$. The manifold $\Mt$ has at least three boundary components, so $\dim \ker i_{*} \geqslant 3$. 
By the relative exact homology sequence 
$\image \d_{*}=\ker i_{*}$, so $H_{2}(\Mt,\d \Mt;\R) \cong \R^{n}$, 
for some $n \geqslant 3$. Let $B$ be the unit ball of $H_{2}(\Mt,\d 
\Mt;\R)$, which is a polyhedron by Theorem \ref{polyhedron}.

Let $L$ be a line with rational slope in $H_{1}(\Tt;\R)$. Then 
$\phi^{-1}(L)$ is an $n-1$ dimensional subspace of $H_{2}(\Mt,\d 
\Mt;\R)$. 

If $\phi^{-1}(L)$ does not intersect the interior of a top 
dimensional 
face of $B$, then neither does any line in $\phi^{-1}(L)$. So any 
point of $\phi^{-1}(L)-\ker \phi$ defines a one dimensional subspace 
$V$, such that $\phi(V)=L$, and $V$ does not pass through the 
interior of a top dimensional face of $B$.

If $\phi^{-1}(L)$ does intersect the interior of a top dimensional 
face of $B$, then the intersection 
has dimension $n-2 \geqslant 1$, so $\phi^{-1}(L)$ must intersect the 
boundary of that face in at least $n-1$ linearly independent points. 
At least one of these points must be non-zero under $\phi$ as $\dim 
\ker \phi = n-2$. This point defines a line $V$ which does not pass 
through the interior a top dimensional face of $B$. We now need to 
show $V$ can be chosen to have rational slope.

We have chosen a preferred basis of $H_{2}(\Mt, \d \Mt;\R)$ in which 
all elements of $H_{2}(\Mt, \d \Mt;\Z)$ are represented by integer 
multiples of the basis elements. The the faces of $B$ are defined by 
linear equations and inequalities with integer coefficients. The 
subspaces $\ker \phi$ and $\phi^{-1}(L)$, are also defined by linear 
equations with integer coefficients, so if the intersection of 
$\phi^{-1}(L)-\ker \phi$ with any face of $B$ is non-empty, it will 
contain a point with rational coefficients. So $V$ can always be 
chosen to have rational slope.
\end{pf}

We can now prove Theorem \ref{main}:

Lemma \ref{onto}  shows that for all slopes $a/b$ on $\Tt$ there is a 
line $V$ in $H_{2}(\Mt,\d \Mt;\R)$ with rational slope such that $\rho 
\d_{*}V=L$, and $V$ does not pass through the interior of a top 
dimensional face of the unit ball. So there is an embedded 
non-separating surface $S$ in $\Mt$, which is not the fiber of a 
fibration, and which has the required boundary slope on $\Tt$. Note 
that we can choose $S$ to be a norm-minimising surface in its 
homology class, i.e. $\chi_{-}(S) = \norm{[S]}$, so $S$ will be 
incompressible. 

However, to apply Theorem \ref{bc} we need a \emph{connected} surface 
with these properties. Suppose $S_{1}$ and $S_{2}$ are two surfaces 
such that $\norm{[S_{1}+S_{2}]} = \norm{[S_{1}]} + \norm{[S_{2}]}$. 
Let $[S_{1}] = q_{1} a_{1}$ and $[S_{2}] = q_{2} a_{2}$, where $0 < 
q_{i} \in \Q$, and $\norm{a_{i}} = 1$. Then the intersection of the 
line through $[S_{1}+S_{2}]$ with the surface of the unit ball is 
given by 
\[ 
\frac{q_{1}}{q_{1}+q_{2}}\,a_{1}+\frac{q_{2}}{q_{1}+q_{2}}\,a_{2}\] 
which lies on the straight line connecting $a_{1}$ and $a_{2}$. 
Therefore the line segment between $a_{1}$ and $a_{2}$ must lie in 
the surface of the unit ball by convexity, so $a_{1}$ and $a_{2}$ lie 
in the same top dimensional face of $B$, though not necessarily in 
its interior.

 Suppose $S$ is not connected, with connected components 
$S_{1},\ldots ,S_{k}$. As $S$ is norm-minimising, $\norm{[S_{1}+\dots 
+S_{k}]} = \norm{[S_{1}]}+\dots +\norm{[S_{k}]}$, so all the lines 
defined by the $S_{i}$ pass through the same top dimensional face of 
$B$. As this face is convex, if any line passes through the interior, 
then the sum of the homology classes will define a line passing 
through the interior, so every component of $S$ defines a line which 
does not pass through a top dimensional face of $B$. As $S$ has 
boundary on $\Tt$, so does at least one of its connected components. 
This component is a connected surface with the correct slope on $\Tt$, 
which is not  a fiber of a fibration.

This completes the proof of Theorem \ref{main}.

\subsection{A Remark on Condition 2 of Theorem \ref{main}}

What happens if $\Tt$ is not a degree one cover of $\d M$? Let 
$p:\Tt \rightarrow \d M$, and let $\mu$ and $\lambda$ be slopes of $\Tt$ 
which cover $m$ and $l$ respectively. If $p(\mu)=q_{1}m$ and 
$p(\lambda)=q_{2}l$, then an embedded surface in $\Mt$ with slope $a/b$ 
on $T$ projects down to an immersed surface with slope 
$q_{1}a/q_{2}b$. However this will only be an embedding on the 
boundary if $q_{1}a$ and $q_{2}b$ are coprime. For the applications 
in this paper, we will always be able to choose $\Tt$ to be a one-fold 
cover.

\subsection{Remarks on Condition 3 of Theorem \ref{main}} 
\label{filling}

In order to apply Theorem \ref{main}, we need to find a boundary 
component $\Tt$ of $\Mt$, such that $\rho:\ker i_{*} \rightarrow 
H_{1}(\Tt;\R)$ is onto. In this section we investigate what happens 
when this projection map is not onto. 
 
Let $M(r)$ denote the manifold obtained by Dehn filling $M$ with 
slope $r$. Let $\Mt(r)$ be the manifold obtained by Dehn filling 
$\Mt$ such that the filling curve on each boundary component of $\Mt$ 
covers the filling curve with slope $r$ on $\d M$.

\begin{lem}{} \label{homology}
Suppose $\Mt$ has at least two boundary components. If the projection 
$\rho:\ker i_{*} \rightarrow H_{1}(\Tt;\R)$ is not onto, then $\dim 
H_{2}(\Mt(0)) \geqslant 2$. 
\end{lem}
\begin{pf}
Suppose $\rho$ is not onto.

Let $\intersecto$ be the \emph{intersection form} on $H_{1}(\d 
\Mt;\R)$. The value of $\intersect{\alpha}{\beta}$ is defined to be the algebraic intersection number of 
representatives of the homology classes of $\alpha$ and $\beta$ in general position. The form $\intersecto$ is a skew-symmetric bilinear form on 
$H_{1}(\d \Mt;\R)$, which is  non-singular on $H_{1}(\d \Mt;\R)$, and $\intersecto \equiv 0$ on $\ker i_{*}$.

Note that $l$ is null homologous in $H_{1}(M;\R)$, so there is a 
surface $S$ in $M$ with boundary parallel to $l$. There is a  
pre-image of this surface, $\St$ in $\Mt$ with boundary on $\Tt$. Take the surface $\St$ to be the entire pre-image of $S$, so $\St$ need not be connected, but it does intersect every boundary component of $\Mt$. The boundary of $\St$ is in $\ker i_{*}$, so $\rho([\d \St])$ is a non-zero 
multiple of $l$. Therefore if $\rho$ is not onto, then its image must be one 
dimensional, generated by $l$. Let $\lambda$ be the slope on $\Tt$ which 
covers $l$. Every element of $\ker i_{*} \cap H_{1}(\Tt;\R)$ must be 
represented by slopes on $\Tt$ parallel to $\lambda$, as the map induced by the covering map
$p_{*}:H_{1}(\Tt;\R) \rightarrow H_{1}(\d M;\R)$ is injective. So 
$\intersect{\lambda}{\alpha} = 0$ for all $\alpha \in \ker i_{*}$, as every $\alpha \in \ker i_{*}$ is represented by some surface with boundary parallel to $\lambda$ on $\Tt$, or else with no boundary at all on $\Tt$. The 
form $\intersecto$ is non-singular, and $\ker i_{*}$ is a maximal 
subspace on which it vanishes. So as $\intersect{\lambda}{\alpha} = 0$ 
for all $\alpha \in \ker i_{*}$ this means that $\lambda \in \ker i_{*}$. Therefore there 
is a surface $S'$ in $\Mt$, such that $\d S'$ is a multiple of 
$\lambda$ on $\Tt$, and the homology of $\d S'$ is zero on all other boundary components of $\d \Mt$. So we can choose this surface $S'$ to have boundary only on $\Tt$.

Dehn fill all boundary components of $M$ and $\Mt$ with slope $0$, to 
produce closed manifolds $M(0)$ and $\Mt(0)$, so that $\Mt(0)$ is a 
branched covering of $M(0)$.

Now consider the surfaces formed from $\St$ and $S'$ by capping off 
their boundaries with meridinal discs of the solid tori filling the 
boundary components of $\Mt$. These two surfaces are non-zero 
elements of $H_{2}(\Mt(0);\R)$. The cover of a meridian of $\d M$ on 
any component of $\d \Mt$ is a non-zero homology class in $\Mt(0)$. 
As there are at least two components of $\d \Mt$, there is a cover of 
the meridian which intersects $\St$ but not $S'$, so these two 
surfaces represent linearly independent elements of 
$H_{2}(\Mt(0);\R)$. Therefore $\dim H_{2}(\Mt(0);\R) \geqslant 2$. 
\end{pf}

Lemma \ref{homology} shows that for a boundary torus $\Tt$ in $\d \Mt$, 
then if the map $\rho : \ker i_{*} \rightarrow H_{1}(\Tt;\R)$ is not 
onto, then there is a surface $S'$ with boundary consisting only of 
parallel copies of $\lambda$, the slope on $\Tt$ which covers $l$. 

Linking numbers, $\lk{\alpha}{\beta} \in \Q$, are defined between 
disjoint simple closed curves which represent elements in the torsion 
subgroup of $H_{1}(\Mt(\infty);\Z)$, i.e $[\alpha] = [\beta] = 0$ in 
$H_{1}(\Mt(\infty);\R)$. The manifold $\Mt(\infty)$ is a branched cover 
of $M(\infty)$, with branch set the core $K$ of the filling torus of 
$\d M$. For certain branched covers, all of the components of $\Kt$, 
the cores of the solid tori filling $\d \Mt$, are in the torsion 
subgroup of $H_{1}(\Mt(\infty);\R)$, so linking numbers are defined 
between them. When this happens, the surface $S'$, considered as a surface in $\Mt(\infty)$, is a surface whose boundary is homologous to some non-zero multiple of the core of $\Tt$. The surface $S'$ doesn't intersect any of the other boundary components of $\d \Mt$, so it doesn't intersect any of the cores of the other tori filling them. Therefore the surface $S'$ 
shows that the core of $T$ must have linking number zero with the 
cores of the other filling tori. This gives the following result:

\begin{lem} \label{linking}
If the projection $\rho:\ker i_{*} \rightarrow H_{1}(\Tt;\R)$ is not 
onto, and linking numbers are defined between the cores of the 
filling tori in $\Mt(\infty)$, then the core of the torus filling $\Tt$ 
has linking number zero with each of the cores of the other filling 
tori.
\end{lem}

\section{Hyperbolic Punctured Torus Bundles.}

We need to show that there is a cover of a hyperbolic punctured torus 
bundle satisfying the conditions of Theorem \ref{main}.

In this section, $M$ will denote the punctured torus bundle, $M=F 
\times 
I / (x,0)\sim (f(x),1)$, where $F$ is a punctured torus, and 
$f:F\rightarrow F$ is a homeomorphism of the punctured torus to 
itself. 

\begin{figure}[h]
\begin{center}
\begin{picture}(60,60)
\put(0,0){\framebox(60,60){}}
\put(32,0){\vector(1,0){0}}
\put(0,32){\vector(0,1){0}}
\put(30,30){\circle{15}}
\put(22.5,30){\circle*{2}}
\put(10,23){{$x_{0}$}}
\put(37.5,31){\vector(0,1){1}}
\put(40,30){{$l$}}
\put(2,30){{$y$}}
\put(30,2){{$x$}}
\end{picture}
\caption{A punctured torus} \label{ptorus}
\end{center}
\end{figure}
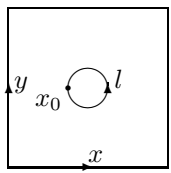

Pick a basepoint $x_{0} \in \d F$. The loops labelled $x$ and $y$ 
give 
a basis for $H_{1}(F;\Z)$, and the loops $l$ and $m=x_{0} \times I$ 
are a  basis for $H_{1}(\d F;\Z)$. This is illustrated in Figure 
\ref{ptorus}.

Each homeomorphism induces an automorphism, $f_{*}$ of $H_{1}(F;\Z) 
\cong \Z \times \Z$. So $f_{*} \in SL_{2}(\Z)$. If $\vert 
\trace(f_{*}) \vert >2$, then the homeomorphism is said to be 
\emph{pseudo-Anosov}, and the resulting manifold is hyperbolic 
\cite{Jorgensen77,Sullivan}. If $\trace(f_{*})=3$, then the manifold is 
the 
figure eight knot exterior and the homology basis can be chosen so 
that the loops $l$ and $m$ are the longitude and 
meridian of the figure eight knot. Note that a punctured torus bundle 
is irreducible and has incompressible boundary. 

We will construct the covering space in the following way: 

The fundamental group of $F$ is the free group on two generators, 
$<x,y \mid>$. Let $\phi:\pi_{1}F \rightarrow \field{3} \oplus 
\field{3}$ be the homomorphism that sends $x \mapsto \col{1}{0}$ and 
$y \mapsto \col{0}{1}$, where $\Z_{3}$ is the cyclic group of order 
$3$. Then $\ker \phi$ defines a $9$-fold regular covering space $\Ft$ 
of $F$. The subgroup $\ker \phi$ is a characteristic subgroup of 
$\pi_{1}F$, so $f_{*}(p_{*}\pi_{1}\Ft)=p_{*}\pi_{1}\Ft$, and there is 
a homeomorphism $\ft:\Ft \rightarrow \Ft$ which covers $f$.

\begin{figure}[h]
\begin{center}
\begin{picture}(80,80)(-20,-20)
\multiput(0,0)(20,0){3}{\multiput(0,0)(0,20){3}{\usebox{\ptorus}}}
\put(-15,5){0}
\put(-15,25){1}
\put(-15,45){2}
\put(10,-15){0}
\put(30,-15){1}
\put(50,-15){2}
\put(30,0){\vector(1,0){0}}
\put(0,30){\vector(0,1){0}}
\put(30,1){\text{$\xt$}}
\put(1,30){\text{$\yt$}}
\end{picture}
\caption{The $3 \times 3$-fold cover} \label{3cover}
\end{center}
\end{figure}
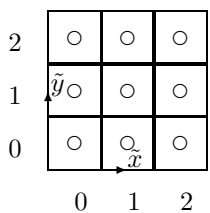

The covering translations of $\Ft$ form a group isomorphic to $\Z_{3} 
\times \Z_{3}$. Choose a basepoint $\xt_{0}$ for $\Ft$ on a boundary 
component such that $p(\xt_{0})=x_{0}$, and label this boundary 
component $L{\col{0}{0}}$. Label each other boundary component by 
the covering translation which maps $L{\col{0}{0}}$ onto it. This 
labeling is illustrated in Figure \ref{3cover}. It is easy to check 
that with respect to the standard homology basis $\ft$ permutes the 
boundary components, as labelled, in the same 
way that the matrix of $f_{*}$ permutes the elements of $\Z_{3} 
\times \Z_{3}$ (The $2 \times 2$ matrix of integers $f_{*}$ acts on 
$\Z_{3} \times \Z_{3}$ in the obvious way). As the monodromy 
homeomorphism $f:F \rightarrow F$ is covered by a 
homeomorphism $\ft:\Ft \rightarrow \Ft$, so $\Mt = \Ft \times I / 
(x,0) \sim (\ft(x),1)$ is a covering space for $M$. Note that $\Mt$ is an irregular covering of $M$ for the punctured torus bundles we are considering, even though $\Ft$ is a regular covering of $F$. 

Note also that as 
$f$ is pseudo-Anosov, $\ft$ will be pseudo-Anosov. 

The top of the cylinder formed by the boundary component 
$L{\col{a}{b}} \times I$ of $\Ft \times I$, is identified with the 
boundary component of $F$ in the position given by 
$f_{*}(\begin{smallmatrix}a \\ b \end{smallmatrix})$. Think of the 
linear map $f_{*}$ acting on the two 
dimensional vector space $\Z_{3} \times \Z_{3}$ as a permutation, 
then the number of boundary components of $\Mt$ is equal to the 
number of cycles of $f_{*}$. Note that as 
$\col{0}{0}$ always gets mapped to itself, the permutation will 
always have at least one  1-cycle. The boundary torus formed by this 
boundary component of $\Ft$ will be chosen to be the boundary torus 
$\Tt$, which is a one-fold covering of $\d M$.

The following table shows the cycles of $f_{*}$. All conjugacy classes of $PSL_{2}(\Z_{3})$ are listed, labelled by a single element chosen from each conjugacy class:  

\begin{center}
\begin{tabular}{rl} 
\\
Conjugacy class of $f_{*}$ & Length of cycles in $\field{3} \times 
\field{3}$ \\
\vspace{3pt}$(\begin{smallmatrix}1 & 0 \\ 0 & 1 \end{smallmatrix})$ & 
1,1,1,1,1,1,1,1,1,1 \\
\vspace{3pt}$(\begin{smallmatrix}1 & 1 \\ 0 & 1 \end{smallmatrix}) 
(\begin{smallmatrix}1 & 2 \\ 0 & 1 \end{smallmatrix})$ & 1,1,1,3,3 \\
\vspace{3pt}$(\begin{smallmatrix}2 & 0 \\ 0 & 2 \end{smallmatrix})$ & 
1,2,2,2,2\\
\vspace{3pt}$(\begin{smallmatrix}2 & 1 \\ 0 & 2 \end{smallmatrix}) 
(\begin{smallmatrix}2 & 2 \\ 0 & 2 \end{smallmatrix})$ & 1,2,6 \\
\vspace{3pt}$(\begin{smallmatrix}0 & 2 \\ 1 & 0 \end{smallmatrix})$ & 
1,4,4 \\
\end{tabular}
\end{center}

So $\Mt$ always has at least three boundary components.

We will use Lemma \ref{homology} to show that the projection 
$\rho:\ker i_{*} \rightarrow H_{1}(\Tt;\R)$ is onto.

Consider Dehn filling all boundary components of $M$ and $\Mt$ with 
slope $0$, to produce closed manifolds $M(0)$ and $\Mt(0)$. For the 
particular covering we have chosen, all pre-images of $l$ cover $l$ 
one-to-one, so the branching index of each Dehn filling is one. 
Therefore the manifold $\Mt(0)$ covers $M(0)$, which is a torus 
bundle with Anosov monodromy, so $\Mt(0)$ must also be a torus bundle 
with Anosov monodromy. In particular, this implies 
$H_{2}(\Mt(0);\R)\cong \R$. However this contradicts Lemma 
\ref{homology}, so in fact $\rho$ must be onto.

So $\Mt$ satisfies the conditions of Theorem \ref{main}. Therefore 
for hyperbolic punctured torus bundles, every slope is a virtually 
embedded boundary slope. This proves Corollary \ref{torus}.

\section{Knots in $S^{2} \times S^{1}$}

We can use Lemma \ref{homology} whenever we know that 
$H_{2}(\Mt(0);\R) \cong \R$. In particular, this happens for all 
finite covers of $S^{2} \times S^{1}$.

Let $K$ be a knot in $S^{2} \times S^{1}$ such that the algebraic 
intersection number of $K$ with $S^{2} \times \{\text{point}\}$, 
$\vert \intersect{K}{S^{2} \times \{\text{pt}\}} \vert = n \geqslant 
3$, and assume further that $M = S^{2} \times S^{1} - K$ is 
hyperbolic. 

Take the $n$-fold cover of $S^{2} \times S^{1}$. The knot $K$ lifts 
to a link with $n$ components, each of which is a degree one cover of 
$K$. As this covering of $M$ doesn't unwrap $\d M$ in the direction 
of the longitude, the covering $p:\Mt \rightarrow M$ extends to a 
covering $\bar p:\Mt(0) \rightarrow M(0)$. 

Suppose the projection $\rho: \ker i_{*} \rightarrow H_{1}(T;\R)$ is 
not onto, then by Lemma \ref{homology}, $\dim H_{2}(\Mt(0);\R) 
\geqslant 2$. But $\Mt(0)$ covers $M(0) = S^{2} \times S^{1}$, so 
$\Mt(0)$ must also be $S^{2} \times S^{1}$. But then 
$H_{2}(\Mt(0);\R) \cong \R$, which gives a contradiction.

Therefore, for these manifolds, every slope is a virtually embedded 
boundary slope. This proves Corollary \ref{s2s1}.

\section{Two-bridge knots}

In this section, the manifold $M$ will always be a hyperbolic knot 
space, i.e. the complement of the interior of a regular neighbourhood 
of the knot $K$ in $S^{3}$. We write $b(\alpha,\beta)$ to denote the 
two bridge knot which gives the lens space $L(\alpha,\beta)$, when 
used as the branch set for a two-fold branched cover of $S^{3}$. If 
$b(\alpha,\beta)$ is a knot, then $\alpha$ is odd. The only closed 
incompressible surfaces in two bridge knots are boundary parallel 
tori \cite{HT}, so all two bridge knots are hyperbolic, except those 
that are torus knots \cite{Thurston2}.

In general, the cores of the filling tori of $\Mt(\infty)$, need not 
be null-homologous in $H_{1}(\Mt;\R)$, so linking numbers need not 
exist between them. However linking numbers  do exist for particular 
classes of branched covers, corresponding to dihedral representations 
of knot groups, which have been extensively studied. For full details 
of all the results used about dihedral coverings see \cite[Chapter 
14]{BZ}.

Suppose we have a representation $\phi:\pi_{1}M \rightarrow D_{2n} = 
\Z_{2} \ltimes \Z_{n}$, with $n$ odd. The fundamental group of $M$ 
can be written as $\Z \ltimes G$, where $\Z$ is generated by a 
meridian of the knot, and $G$ is the commutator subgroup of 
$\pi_{1}M$. Therefore $m$ must get mapped onto a reflection in 
$D_{2n}$. The longitude $l$ is in the commutator subgroup of 
$\pi_{1}M$, so $\phi(l)$ is in the $\Z_{n}$ subgroup of rotations. 
But $l$ and $m$ commute, so $\phi(l)=1$. So there is a regular 
$2n$-fold covering space $M_{2n}$ corresponding to $\ker \phi$, which 
has $n$ boundary components, each of which is a two-fold cover of $\d 
M$. Let $A$ be the $\Z_{2}$ subgroup generated by $\phi(m)$. Then 
$\phi^{-1}(A)$ generates an irregular $n$-fold covering $M_{n}$ of 
$M$. 

We need to know how many boundary components the cover $M_{n}$ has. Choose a point $x$ in $\d M$. The group $\pi_{1}M$ acts on $p^{-1}(x)$ on the right by path lifting, i.e. if $\xt \in p^{-1}(x)$, and $\alpha$ is a loop based at $x$ in $M$, then $\xt \alpha = \alphat(1)$, where $\alphat$ is the unique lift of $\alpha$ such that $\alphat(0) = \xt$. As a right $\pi_{1}M$ space, $p^{-1}(x)$ is isomorphic to the space of right cosets of $p_{*}\pi_{1}M_{n}$ in $\pi_{1}M$. This in turn is isomorphic to the space of right cosets of $A$ in $D_{2n}$, as an element $[\alpha] \in \pi_{1}M$ acts on $A$ by right multiplication by $\phi([\alpha])$. Choose two elements $x_{1},x_{2}$ of $p^{-1}(x)$ and label them by the right $A$ cosets, $Aa_{1},Aa_{2}$, they correspond to. The two elements $x_{1},x_{2}$ of $p^{-1}(x)$ lie in the same boundary component of $\d M_{n}$, if and only if there is a path $\alpha$ in $\d M_{n}$ connecting them. This path projects down to a loop $p(\alpha)$ in $\d M$, which gets mapped into $A$ by $\phi$, as $\phi(\pi_{1}\d M) = A$. So $Aa_{1}\phi([p(\alpha)]) = Aa_{2}$. This element $\phi([p(\alpha)])$ of $A$ must map one coset to the other, and if any element of $A$ does so, then there is a corresponding path in $\d M_{n}$ connecting the two points. Therefore two elements of $p^{-1}(x)$ lie in the same boundary component if and only if their corresponding cosets lie in the same $(A,A)$-double coset of $D_{2n}$. Furthermore the order of the covering of each boundary component is given by the number of cosets of $A$ in each double coset.

A simple calculation shows that the number of $(A,A)$-double 
cosets in $D_{2n}$ is $(n+1)/2$, all of which contain two cosets of 
$A$, except $A$ itself. Therefore the cover $M_{n}$ has $(n+1)/2$ 
boundary components, one of which is a one-fold cover. Choose this 
boundary component to be $\Tt$, which is covered two-to-one by a single 
boundary component of $M_{2n}$. There are $n$ cosets of $A$ in 
$D_{2n}$, so $M_{2n}$ has $n$ boundary components. Therefore all the 
boundary components of $M_{n}$ except $\Tt$ are covered by two boundary 
components of $M_{2n}$.

\begin{thm}{\cite[Theorem 14.8]{BZ}}
There is a surjective homomorphism $\phi:\pi_{1}M \rightarrow 
D_{2p}$, if and only if the prime $p$ divides the order of 
$H_{1}(C_{2})$, where $C_{2}$ is the two-fold branched cover of 
$S^{3}$, branched over the knot $K$. If $p$ does not divide the 
second torsion coefficient of $H_{1}(C_{2})$, then all such 
representations are equivalent.
\end{thm}

\begin{thm}{\cite[Proposition 14.16]{BZ}}
If there is exactly one class of equivalent dihedral homomorphisms 
$\phi:\pi_{1}M \rightarrow D_{2p}$, then linking numbers are defined 
between the core curves of $M_{2p}$, and also $M_{p}$. 
\end{thm}

These linking numbers, when they exist, have been computed for many 
of the knots in the knot tables.  However in the case of two bridge 
knots, they always exist.

\begin{thm}{\cite{Burde}} \label{burde}
Let $M$ be the knot space of the two bridge knot $b(\alpha,\beta)$. 
There is a dihedral representation $\phi:\pi_{1}M \rightarrow 
D_{2\alpha}$. Linking numbers are defined in the branched covers 
corresponding to $\ker \phi$, and $\phi^{-1}(A)$, for any reflection 
subgroup $A \cong \Z_{2}$, even if $\alpha$ is not prime. The linking 
numbers in $M_{2\alpha}$ are $\pm 1$ for all pairs of core curves in 
$M_{2\alpha}$.
\end{thm}

The covers $M_{2\alpha}$ and $M_{\alpha}$ both give rise to branched 
covers of $S^{3}$, namely $M_{2\alpha}(\infty)$ and 
$M_{\alpha}(\infty)$, in the notation of Section \ref{filling}. Note 
that $M_{2\alpha}(\infty)$ is a two-fold branched cover of 
$M_{\alpha}(\infty)$, branched over $t$, the core of $\Tt$. 
Let $c$ be some other core curve in $M_{\alpha}(\infty)$, which will 
have two pre-images $\tilde c_{1}$ and $\tilde c_{2}$ in 
$M_{2\alpha}$. If $\lk{t}{c} = 0$, then $\lk{\tilde t}{\tilde c_{i}}$ 
will also be zero for each $i$,
but by Theorem \ref{burde}, $\lk{\tilde t}{\tilde c_{i}} = \pm 1$ for 
both pre-images of $c$. Therefore $\lk{t}{c} \not = 0$ for all core 
components $c \not = t$, so by Lemma \ref{linking}, $\rho:\ker i_{*} 
\rightarrow H_{1}(\Tt;\R)$ is onto. If $p \geqslant 5$, then 
$M_{\alpha}$ has at least three boundary components, so by Theorem 
\ref{main}, every boundary slope is a virtually embedded boundary 
slope. The only two-bridge knot with $\alpha < 5$ is the trefoil, 
which is not hyperbolic, so this proves Corollary \ref{twobridge}.

Table III in \cite{BZ} lists linking invariants of knots, from which 
the linking numbers can easily be computed. This shows that many of 
the hyperbolic knots in the tables have virtually immersed boundary 
slopes of every slope.

\section{Acknowledgments}
I would like to thank my PhD advisor Daryl Cooper for his continuing 
guidance and support.

\noindent \\ 
J. Maher:\\
Department of Mathematics\\
UCSB CA 93106\\
USA\\
email: maher@math.ucsb.edu \\
fax (via UCSB Math Dept.): (805) 893 2385

\end{document}